\documentclass[11pt, twoside]{article}
\usepackage{amssymb}
\usepackage{amsmath, amssymb}\usepackage{cite}\usepackage{mathrsfs}\usepackage[amsmath, thmmarks]{ntheorem}
\usepackage{listings}
\usepackage[titletoc]{appendix}\allowdisplaybreaks
\usepackage{multirow}

\makeatletter
\renewcommand\theequation{\thesection.\@arabic\c@equation}
\makeatother

\newtheorem{thm}{Theorem}[section]%
\newtheorem{lem}[thm]{Lemma}%
\newtheorem{Con}[thm]{Conjecture}%
\input cyracc.def

\def\f{\noindent}

\def\demo{\f{\bf Proof}\hskip10pt}

\def\qed{\hfill $\Box$}

\begin{document}

\title{{\bf A note on the Cuntz algebra automorphisms}}
\footnotetext{The work was supported by 2023 Excellent Science and Technology Innovation Team of Jiangsu Province Universities (Real-time Industrial Internet of Things).\\
E-mail addresses: Junyao$_{-}$Pan@126.com}

\author{{\bf Junyao Pan}\\
{\footnotesize School of Cyber Science and Engineering, Wuxi University, Wuxi, Jiangsu, 214105, P. R. China}
}

\date{}
\maketitle

%
%
%
%

\noindent{\small {\bf Abstract:} Permutative automorphisms of the Cuntz algebras $\mathcal{O}_n$ are in bijection with the stable permutations of $[n]^k$. They are also the elements of the restricted Weyl group of $Aut(\mathcal{O}_n)$. In this note, we characterize a class of stable involutions of $[n]^2$. More precisely, we prove Conjecture 12.2 of Brenti and Conti [Adv. Math. 381 (2021), p. 60], and thus providing a new family (with $6$ degrees of freedom) of automorphisms of the Cuntz algebras $\mathcal{O}_n$ for any $n>1$.

\vskip0.2cm
\noindent{\small {\bf Keywords:} Cuntz algebra; Permutative automorphism; Stable permutation.


\section {Introduction}

Throughout this paper, $[n]:=\{1,2,...,n\}$; $[n]^{t}$ stands for the cartesian product of $t$ copies of $[n]$, where $t \in N$; $S([n]^{t})$ expresses the symmetric group on $[n]^{t}$.

A $C^*$-algebra is a norm closed self-adjoint sub-algebra of the bounded operators on a Hilbert space. Due to a theorem of Gelfand, Naimark and Segal, one can alternatively describe $C^*$-algebras axiomatically as complex Banach algebras with an involution that satisfies $\|a^*a\|=\|a\|^2$. Actually, $C^*$-algebras were first introduced for providing a suitable environment for a rigorous approach to quantum theories \cite{HAA}, and more recently have been related to such diverse fields as operator theory, group representations, topology, quantum mechanics, non-commutative geometry and dynamical systems \cite{CON}. On the other hand, the symmetries of $C^*$-algebras are provided by automorphisms in the most classical sense. Despite the fact that many general results and constructions about automorphisms of $C^*$-algebras have been obtained in \cite{FA,PE} and that they are used extensively, not very much is known about the automorphism group of most $C^*$-algebras. In addition, many other interesting (simple) examples of $C^*$-algebras emerged in the 1970s and 1980s, such as Cuntz \cite{C1} invented his Cuntz algebras $\mathcal{O}_n$. These are, for $2 \leq n < \infty$, generated by $n$ isometries $s_1, s_2,..., s_n$ satisfying the relation $I = s_1s^*_1+s_2s^*_2+\cdot\cdot\cdot+s_ns^*_n$. For $n = \infty$, this relation is replaced with the relation that the support projections $s_js^*_j$ are mutually orthogonal. Actually, they have been attracting much attention since their appearance in the seminal paper, such as \cite{BJ,C2,C3,CK,DR1,DR2}. Moreover, the study of automorphisms of the Cuntz algebras soon appeared quite intriguing and revealed many interesting facets. Notably, following the deep insight by Cuntz \cite{C4}, in a series of papers \cite{CHS1,CHS2,CS,CKS} a general theory of restricted Weyl groups for $Aut(\mathcal{O}_n)$ has been studied both from a theoretical viewpoint as well as from the perspective of constructing explicit examples. In order to explain the connection with the main topic of this note, we next introduce some more terminology.

 Let $\mathcal{O}_n$ be the Cuntz algebra with $n \geq2$, $\mathcal{F}_n \subset\mathcal{O}_n$ the so-called core-UHF $C^*$-subalgebra generated by a nested family of subalgebras $\mathcal{F}^k_n$ that is isomorphic to the algebra of complex matrices $M_{n^k}$ and $\mathcal{D}_n$ be the $C^*$-subalgebra generated by the family of subalgebras $\mathcal{D}^k_n$ of $\mathcal{F}^k_n$ that is isomorphic to the algebra of diagonal matrices in $M_{n^k}$. Then, following the insight by Cuntz, the reduced Weyl group of $Aut(\mathcal{O}_n)$ is defined as $$Aut(\mathcal{O}_n,\mathcal{F}_n)\cap Aut(\mathcal{O}_n,\mathcal{D}_n)/Aut_{\mathcal{D}_n}(\mathcal{O}_n),$$ where $Aut(\mathcal{O}_n, X)$ is the subgroup of $Aut(\mathcal{O}_n)$ consisting of the automorphisms which leave $X$ invariant, while $Aut_X(\mathcal{O}_n)$ is the one of those which fix $X$ pointwise. It is well-known that unital $*$-endomorphisms of $\mathcal{O}_n$ are in bijective correspondence with unitaries in $\mathcal{O}_n$, call this bijection $u\mapsto\lambda_u$. With some more work the reduced Weyl group can be further identified with the set of automorphisms $\lambda_u$ of $\mathcal{O}_n$ induced by the so-called permutative unitaries $u \in \cup_{k\geq1}\mathcal{F}^k_n$. Moreover, for general permutative unitaries $u \in \mathcal{F}^k_n$, for any $k$, it was shown in \cite{CHS2}, that $\lambda_u$ is an automorphism precisely when the sequence of unitaries $(\varphi^r(u^*)\cdot\cdot\cdot\varphi(u^*)u^*\varphi(u)\cdot\cdot\cdot\varphi^r(u))_{r\geq0}$ in $\mathcal{F}_n$ eventually stabilizes, where the endomorphism $\varphi$ of $\mathcal{F}_n$ corresponds, in the isomorphism between $M_{n^k}$ and $M_n\otimes M_n\otimes\cdot\cdot\cdot\otimes M_n$ with $k$ factors, to the tensor shift map $x \mapsto 1_{M_n}\otimes x$. By coherently identifying permutative unitaries in $\mathcal{F}^k_n$ with permutation matrices in $M_{n^k}$ and thus with permutations of the set $\{1, ..., n^k\}$, and finally with permutations of the set $[n]^k$ by lexicographic ordering, one is lead to the class of \emph{stable permutations} that were defined by Brenti and Conti \cite{BC}. Notably, Brenti and Conti \cite{BC} pointed out that $\lambda_u$ is an automorphism of $\mathcal{O}_n$ if and only if $u$ is stable, in particular, the reduced Weyl group of $\mathcal{O}_n$ is $\{\lambda_u : u \in S([n]^r), r \in N, u~{\rm{ stable}}\}$. This opened up the investigations of the reduced Weyl groups of the Cuntz algebras $\mathcal{O}_n$ from a combinatorial point of view. Following the deep insight, Brenti and Conti and Nenashev continued to deep these investigations and focus on the explicit construction of restricted Weyl group elements using combinatorial techniques in \cite{BCN1,BCN2}. Next we recall some notions and notations about stable permutations, for more detailed information see \cite{BC}.

 Given two permutations $u\in S([n]^t)$ and $v\in S([n]^r)$ with $n\in N$ and $n\geq2$, the \emph{tensor product} of $u$ and $v$ is the permutation $u\otimes v\in S([n]^{t+r})$ defined by $$(u\otimes v)(\alpha,\beta):=(u(\alpha),v(\beta))$$ for all $\alpha\in[n]^t$ and $\beta\in[n]^r$. For a permutation $u \in S([n]^t)$, define a sequence of permutations $\Psi_k(u) \in S([n]^{t+k})$, $k \geq 0$ by setting $\Psi_0(u) := u^{-1}$ and, for $k \in N$, $$\Psi_k(u)=\prod_{i = 0}^{k}\underbrace{1\otimes\cdot\cdot\cdot\otimes1}_{k-i}\otimes u^{-1}\otimes\underbrace{1\otimes\cdot\cdot\cdot\otimes1}_{i}\prod_{i = 1}^{k}\underbrace{1\otimes\cdot\cdot\cdot\otimes1}_{i}\otimes u\otimes\underbrace{1\otimes\cdot\cdot\cdot\otimes1}_{k-i},$$ where $1$ denotes the identity of $S_n := S([n])$. Then $u$ as above is said to be \emph{stable} if there exists some integer $k \geq 1$ such that $$\Psi_{k+l}(u)=\Psi_{k-1}(u)\otimes\underbrace{1\otimes\cdot\cdot\cdot\otimes1}_{l+1}~{\rm{for ~all}} ~l\geq0,$$
and the least such value of $k$ is called the \emph{rank} of $u$, denoted by $rank(u)$. One easily checks that if $t=1$ then all permutations in $S_n$ are stable of rank $1$. However, it seems to be extremely difficult to determine whether $u$ is stable for $t>1$. Brenti and Conti \cite{BC} posed a number of interesting open problems, including three conjectures, two of which have now been solved, and the unresolved conjecture is as follows:

\begin{Con}\label{pan1-1}\normalfont([2, Conjecture 12.2])
Let $(a_1, b_1),(a_1, b_2),(a_2, b_3),(a_2, b_4) \in [n]^2$ be distinct, $a_1 \neq a_2$. Then
$$u := ( (a_1, b_1),(a_1, b_2)) ( (a_2, b_3),(a_2, b_4))$$
is stable of rank $1$ if and only if one of the following conditions is satisfied:\\\\
$(i)$ $\{a_1, a_2\} \cap \{b_1, b_2, b_3, b_4\}= \emptyset$;\\
$(ii)$ $\{a_1, a_2\} = \{b_1, b_2, b_3, b_4\}$.
\end{Con}

Brenti and Conti \cite{BC} pointed out that the Conjecture\ \ref{pan1-1} holds if one can show that either $(i)$ or $(ii)$ is a necessary condition for $u$ to be stable of rank $1$. In this note, we give a simple proof by contradiction to show that either $(i)$ or $(ii)$ is a necessary condition for $u$ to be stable of rank $1$. So we obtain the following theorem, thus providing a new family (with $6$ degrees of freedom) of automorphisms of the Cuntz algebras $\mathcal{O}_n$ for any $n>1$.

\begin{thm}\label{pan1-2}\normalfont
Let $(a_1, b_1),(a_1, b_2),(a_2, b_3),(a_2, b_4)$ be four distinct elements in $[n]^2$ with $a_1 \neq a_2$. Then the involution $((a_1, b_1),(a_1, b_2)) ((a_2, b_3),(a_2, b_4))\in S([n]^2)$ is stable of rank $1$ if and only if one of the following conditions is satisfied:\\\\
$(i)$ $\{a_1, a_2\} \cap \{b_1, b_2, b_3, b_4\}= \emptyset$;\\
$(ii)$ $\{a_1, a_2\} = \{b_1, b_2, b_3, b_4\}$.
\end{thm}

\section {Proof of Theorem 1.2}
Recall some notions and notations. If $\pi\in S_n$ with $\pi(i_1)=i_2,\pi(i_2)=i_3,...,\pi(i_{r-1})=i_r,\pi(i_r)=i_1$ and $\pi(i)=i$ for all $i\in[n]\setminus\{i_1,i_2,...,i_r\}$, then $\pi$ is called a \emph{$r$-cycle}, and is denoted by $\pi=(i_1,i_2,...,i_r)$. In particular, if $r=2$ then $\pi$ is called a \emph{transposition}. If $\tau=(j_1,j_2,...,j_s)$ and $\{i_1,i_2,...,i_r\}\cap\{j_1,j_2,...,j_s\}=\emptyset$, then $\pi$ and $\tau$ are called two \emph{disjoint} cycles. It is well-known that every permutation can be written as a product of some disjoint cycles. If a permutation is a product of some disjoint transpositions, then this permutation is called an \emph{involution}. So $u$ in Conjecture\ \ref{pan1-1} is an involution. Next we provide a key lemma in our proof.

\begin{lem}\label{pan2-0}\normalfont([2, Proposition 4.5])
Let $v \in S([n]^2)$. Then $v$ is stable of rank $1$ if and only if $v$ satisfies the equation $(v \otimes 1)(1 \otimes v) = (1 \otimes v)(v \otimes 1)$ in $S([n]^3)$.
\end{lem}
Finally, we state some facts that will be repeatedly applied. Since $(a_1, b_1),(a_1, b_2), (a_2, b_3),(a_2, b_4)$ are four distinct elements, it follows that $b_1\neq b_2$ and $b_3\neq b_4$. Additionally, $((a_1, b_1),(a_1, b_2))$ and $((a_2, b_3),(a_2, b_4))$ are symmetrical, and so only one needs to be considered. Similarly, $(a_1, b_1)$ and $(a_1, b_2)$, $(a_2, b_3)$ and $(a_2, b_4)$ are also symmetrical, respectively.

Now we present the main idea of the proof. We prove the necessary conditions for Conjecture\ \ref{pan1-1} by using the method of contradiction. Thus, we assume that $\{a_1, a_2\} \cap \{b_1, b_2, b_3, b_4\}\neq \emptyset$ and $\{a_1, a_2\} \neq \{b_1, b_2, b_3, b_4\}$. Note that there exist two cases, that is, $|\{a_1, a_2\} \cap \{b_1, b_2, b_3, b_4\}|=1$ and $\{a_1, a_2\}\subset \{b_1, b_2, b_3, b_4\}$. However, for each case, we construct an $\alpha\in[n]^3$ such that $$(u \otimes 1)(1 \otimes u)(\alpha)\neq (1 \otimes u)(u \otimes 1)(\alpha),$$ this is contradict to Lemma\ \ref{pan2-0}. Hence, our assumption is not valid. So we arrive at Theorem\ \ref{pan1-2}. Next, we start to discuss the necessary conditions for Conjecture\ \ref{pan1-1} in these two situations.

\begin{lem}\label{pan2-1}\normalfont
Suppose that $(a_1, b_1),(a_1, b_2),(a_2, b_3),(a_2, b_4)$ are four distinct elements in $[n]^2$ with $a_1 \neq a_2$ and $|\{a_1, a_2\} \cap \{b_1, b_2, b_3, b_4\}|=1$. Let $$u:=( (a_1, b_1),(a_1, b_2)) ( (a_2, b_3),(a_2, b_4))\in S([n]^2).$$ Then $u$ is not stable of rank $1$.
\end{lem}
\demo We see that either $\{a_1, a_2\} \cap \{b_1, b_2, b_3, b_4\}=\{a_1\}$ or $\{a_1, a_2\} \cap \{b_1, b_2, b_3, b_4\}=\{a_2\}$. According to the symmetry, it suffices to consider $\{a_1, a_2\} \cap \{b_1, b_2, b_3, b_4\}=\{a_1\}$. Moreover, we note that there are three cases, those are: $(1)$ $a_1\in\{b_1,b_2\}$ and $a_1\not\in\{b_3,b_4\}$; $(2)$ $a_1\not\in\{b_1,b_2\}$ and $a_1\in\{b_3,b_4\}$; $(3)$ $a_1\in\{b_1,b_2\}$ and $a_1\in\{b_3,b_4\}$. We claim that for each case, there exists an $\alpha\in[n]^3$ such that $(u \otimes 1)(1 \otimes u)(\alpha)\neq (1 \otimes u)(u \otimes 1)(\alpha)$.

Case $(1)$: $a_1\in\{b_1,b_2\}$ and $a_1\not\in\{b_3,b_4\}$. Note that either $b_1=a_1$ or $b_2=a_1$. By symmetry, it suffices to consider the case that $b_1=a_1$. In this case, we pick $\alpha=(a_1,b_1,b_2)$. By definition of tensor product, we see that $(u \otimes 1)(\alpha)=(a_1,b_2,b_2)$. Since $b_1\neq b_2$, we have $b_2\neq a_1$. In addition, $\{a_1, a_2\} \cap \{b_1, b_2, b_3, b_4\}=\{a_1\}$ indicates $b_2\neq a_2$, and thus $(1 \otimes u)((a_1,b_2,b_2))=(a_1,b_2,b_2)$. Thereby, $$(1 \otimes u)(u \otimes 1)(\alpha)=(a_1,b_2,b_2).$$
Since $b_1=a_1$, we deduce that $(1 \otimes u)(\alpha)=(a_1,a_1,b_1)$ and $(u \otimes 1)((a_1,a_1,b_1))=(a_1,b_2,b_1)$, and therefore, $$(u \otimes 1)(1 \otimes u)(\alpha)=(a_1,b_2,b_1).$$
It follows from $b_1\neq b_2$ that $(u \otimes 1)(1 \otimes u)(\alpha)\neq (1 \otimes u)(u \otimes 1)(\alpha)$.

Case $(2)$: $a_1\not\in\{b_1,b_2\}$ and $a_1\in\{b_3,b_4\}$. Similarly, it suffices to consider the case that $b_3=a_1$. Pick $\alpha=(a_2,a_1,b_1)$. Note that $b_3=a_1$, $b_4\neq a_1$ and $b_4\neq a_2$. Then by definition of tensor product, we deduce that $$(1 \otimes u)(u \otimes 1)(\alpha)=(1 \otimes u)((a_2,b_4,b_1))=(a_2,b_4,b_1)$$ and $$(u \otimes 1)(1 \otimes u)(\alpha)=(u \otimes 1)((a_2,a_1,b_2))=(a_2,b_4,b_2).$$ Clearly, $(u \otimes 1)(1 \otimes u)(\alpha)\neq (1 \otimes u)(u \otimes 1)(\alpha)$ as $b_1\neq b_2$.

Case $(3)$: $a_1\in\{b_1,b_2\}$ and $a_1\in\{b_3,b_4\}$. By symmetry, we may assume that $a_1=b_1=b_3$. Take $\alpha=(a_1,b_1,b_2)$. Similarly, it follows from $b_2\neq a_1$, $b_2\neq a_2$ and $a_1=b_1$ that $$(1 \otimes u)(u \otimes 1)(\alpha)=(1 \otimes u)((a_1,b_2,b_2))=(a_1,b_2,b_2)$$ and $$(u \otimes 1)(1 \otimes u)(\alpha)=(u \otimes 1)((a_1,a_1,b_1))=(a_1,b_2,b_1).$$ Obviously, $(u \otimes 1)(1 \otimes u)(\alpha)\neq (1 \otimes u)(u \otimes 1)(\alpha)$ due to $b_1\neq b_2$.

According to the above discussions, we deduce this lemma.  \qed

\begin{lem}\label{pan2-2}\normalfont
Suppose that $(a_1, b_1),(a_1, b_2),(a_2, b_3),(a_2, b_4)$ are four distinct elements in $[n]^2$ with $a_1 \neq a_2$ and $\{a_1, a_2\}\subset \{b_1, b_2, b_3, b_4\}$. Let $$u:=( (a_1, b_1),(a_1, b_2)) ( (a_2, b_3),(a_2, b_4))\in S([n]^2).$$ Then $u$ is not stable of rank $1$.
\end{lem}
\demo Note that there exist four cases, those are: $(1)$ $\{a_1,a_2\}=\{b_1,b_2\}$ and $\{a_1,a_2\}\neq\{b_3,b_4\}$; $(2)$ $\{a_1,a_2\}\neq\{b_1,b_2\}$ and $\{a_1,a_2\}=\{b_3,b_4\}$; $(3)$ $\{a_1,a_2\}\cap\{b_1,b_2\}=\{a_1\}$ and $\{a_1,a_2\}\cap\{b_3,b_4\}=\{a_2\}$; $(4)$ $\{a_1,a_2\}\cap\{b_1,b_2\}=\{a_2\}$ and $\{a_1,a_2\}\cap\{b_3,b_4\}=\{a_1\}$. Apply symmetry, it suffices to consider one of $(1)$ and $(2)$. We claim that for each case of $(1)$ and $(3)$ and $(4)$, there exists an $\alpha\in[n]^3$ such that $(u \otimes 1)(1 \otimes u)(\alpha)\neq (1 \otimes u)(u \otimes 1)(\alpha)$.

Case $(1)$: $\{a_1,a_2\}=\{b_1,b_2\}$ and $\{a_1,a_2\}\neq\{b_3,b_4\}$. Since $\{a_1,a_2\}=\{b_1,b_2\}$, we note that either $a_1=b_1$, $a_2=b_2$ or $a_1=b_2$, $a_2=b_1$. However, $((a_1,a_1),(a_1,a_2))=((a_1,a_2),(a_1,a_1))$, and so we can assume that $a_1=b_1$ and $a_2=b_2$. In addition, either $b_3\notin\{a_1,a_2\}$ or $b_4\notin\{a_1,a_2\}$. By symmetry, it suffices to consider $b_4\notin\{a_1,a_2\}$. Pick $\alpha=(a_1,b_2,b_3)$. Since $a_2=b_2$, it follows that $(1 \otimes u)(\alpha)=(a_1,a_2,b_4)$ and $(u \otimes 1)((a_1,a_2,b_4))=(a_1,b_1,b_4)$, and thus $$(u \otimes 1)(1 \otimes u)(\alpha)=(u \otimes 1)((a_1,a_2,b_4))=(a_1,b_1,b_4).$$
One easily checks that
\begin{align*}
(1 \otimes u)(u \otimes 1)(\alpha)=(1 \otimes u)((a_1,b_1,b_3))&=
  \begin{cases}
    (a_1,a_1,b_3) & \text{if } b_3 \notin\{a_1,a_2\}, \\
    (a_1,a_1,b_2) & \text{if } b_3 =a_1, \\
    (a_1,a_1,b_1) & \text{if } b_3 =a_2.
  \end{cases}
\end{align*}
Since $b_4\neq b_3$ and $b_4\notin\{a_1,a_2\}$, we have $(u \otimes 1)(1 \otimes u)(\alpha)\neq (1 \otimes u)(u \otimes 1)(\alpha)$.

Case $(3)$: $\{a_1,a_2\}\cap\{b_1,b_2\}=\{a_1\}$ and $\{a_1,a_2\}\cap\{b_3,b_4\}=\{a_2\}$. By symmetry, we assume that $a_1=b_1$ and $a_2=b_3$. Pick $\alpha=(a_1,b_1,b_2)$. Note $a_1\neq b_2$ and $a_2\neq b_2$ because $b_1\neq b_2$ and $\{a_1,a_2\}\cap\{b_1,b_2\}=\{a_1\}$. Then by definition of tensor product, we deduce that $$(1 \otimes u)(u \otimes 1)(\alpha)=(1 \otimes u)((a_1,b_2,b_2))=(a_1,b_2,b_2)$$
and $$(u \otimes 1)(1 \otimes u)(\alpha)=(u \otimes 1)((a_1,a_1,b_1))=(a_1,b_2,b_1).$$ By $b_1\neq b_2$, we see that $(u \otimes 1)(1 \otimes u)(\alpha)\neq (1 \otimes u)(u \otimes 1)(\alpha)$.

Case $(4)$: $\{a_1,a_2\}\cap\{b_1,b_2\}=\{a_2\}$ and $\{a_1,a_2\}\cap\{b_3,b_4\}=\{a_1\}$. Apply symmetry, we may assume that $a_1=b_3$ and $a_2=b_1$. Pick $\alpha=(a_1,b_1,b_3)$. Note $b_2\neq a_2$ as $a_2=b_1$ and $b_1\neq b_2$. One easily checks that $$(1 \otimes u)(u \otimes 1)(\alpha)=(1 \otimes u)((a_1,b_2,b_3))=(a_1,b_2,b_3)$$ and $$(u \otimes 1)(1 \otimes u)(\alpha)=(u \otimes 1)((a_1,a_2,b_4))=(a_1,b_2,b_4).$$ It follows from $b_3\neq b_4$ that $(u \otimes 1)(1 \otimes u)(\alpha)\neq (1 \otimes u)(u \otimes 1)(\alpha)$.

According to the above discussions, we deduce this lemma.  \qed

Up to now we have completed the proof of Theorem\ \ref{pan1-2}.

\section{Acknowledgement}

We are very grateful to the anonymous referees for their useful suggestions and comments.

\end{document}